\documentclass{article}
\usepackage{amsmath}
\usepackage{amssymb}
\usepackage{amsfonts}
\newcommand{\bq}{\mathbf{Q}}
\newcommand{\bz}{\mathbf{Z}}

\newcommand{\bd}{\begin{definition}}
\newcommand{\ed}{\end{definition}}
\usepackage{commath} 
\usepackage{csquotes} 
\title{An Infinite, Two-parameter Family of Polynomials with Factorization 
Similar to  $X^m-1$}

\author{M Krithika and P Vanchinathan*}
\begin{document}
\maketitle
\begin{center}
	Division of Mathematics\\
	VIT University\\
	Vandalur--Kelambakkam Road\\
	Chennai, 600 127 INDIA\\
\end{center}
\begin{abstract}
	For a suitable irreducible \textit{base}  polynomial $f(x)\in \bz[x]$ of degree $k$, a family of polynomials  $F_m(x)$ depending on $f(x)$ is constructed with the properties:
	\begin{itemize}
		\item[(i)] there is exactly one factor $\Phi_{d,f}(x)$ for $F_m(x)$ for each divisor $d$ of $m$;
		\item[(ii)] $\deg(\Phi_{d,f}(x))=\varphi(d)\cdot\deg(f)$ generalizing the factorization of $x^m-1$ into cyclotomic polynomials;  
		\item[(iii)] when the base polynomial $f(x) = x-1$ this $F_m(x)$ coincides with $x^m-1$.
	\end{itemize}
As an application, irreducible polynomials of degree 12, 24, 24 are 
constructed having Galois groups of order matching their degrees and isomorphic to 
	$S_3 \oplus C_2 , 
S_3 \oplus C_2\oplus C_2$ 
	and  $S_3 \oplus C_4$ respectively.
\end{abstract}

\section{Introduction}
Cyclotomic polynomials occupy a central role in algebraic number theory:
any single root of them generates a Galois extension of the field of rational
 numbers 
with abelian Galois group. And conversely 
the celebrated result of Kronecker and Weber states that 
any abelian extension of $\bq$ is contained in such 
extensions. The product of these polynomials $\Phi_d(x)$ over  $d$ chosen as divisors of
a pre-determined $m$ gives irreducible factorization of $x^m-1$.
In this paper we arrive at one generalization of such a factorization 
 which leads us to new irreducible polynomials that may be 
regarded as some kind of relativised cyclotomic polynomials (relative to
a given fixed irreducible polynomial).

There have been many attempts to find generalizations earlier. One by 
K.\ Nageswara Rao 
\cite{knr} arises by generalizing the notion of
the set of divisors of a number and using the corresponding M\"obius inversion.

C.\ Kimberling \cite{kimberling}  has generalized it to a  polynomial in three variables
and relates them to hyperbolic sine functions, as opposed to classical cyclotomic polynomials
where the  roots are expressible in terms of trigonometric functions.

In Kwon, Lee and Lee \cite{kll} the authors deal with a natural choice of primitive
elements of subfields of cyclotomic fields, and  their minimal polynomials are
presented as a generalization.

In this paper our  approach is to relativise  the concept which  has 
spawned  new kinds of polynomials different from all the
above.

\section{Preliminaries}
Conventions: 
We take the field of rational numbers, $\bq$ as our base field.
All the algebraic numbers  may  be taken as algebraic 
integers with no loss of generality. 
When we talk of degree or conjugates of an algebraic number it is all with respect to $\mathbf{Q}$, 
unless stated otherwise. 

We use the standard notation $\Phi_m(x)$ for the classical
$m$-th cyclotomic polynomial whose  degree is  given by 
the Euler totient function $\varphi(m)$. 
We write $\mathrm{Gal\,} \big(f (x)\big)$ for the Galois group of 
a polynomial $f(x)$ (which is tacitly assumed to be over $\bq$).

First we define the notion of uniform degree which is very handy for our purposes.

\noindent
\textbf{Uniform Degree}:\quad
Let $\alpha$ be an algebraic number.  We say $\alpha$ is of uniform degree, 
if for every  integer $m\neq0$, the degree of  $\alpha^m$ as an algebraic 
number is the 
same as the degree of $\alpha$.

\noindent
\textbf{Example}:\quad For any square-free integer $d$, 
the algebraic number $1+\sqrt d$ is of uniform degree $2$.

Besides the above example there are many more.
The following result assures us such numbers are aplenty in any number field.

\vskip6pt
\noindent
\textbf{Lemma 1.}\quad \textit{For any algebraic number $\alpha$, say of degree $d$,
there is an integer  $k\in \mathbf{Z}$ such that   $\alpha +k$ has uniform degree $d$}.

\vskip6pt
Proof:\quad Let $\alpha=\alpha_1,\alpha_2,\ldots,\alpha_d$ be all the conjugates 
of $\alpha$. 

Suppose that $\alpha^n$ is ``defective'' in the sense that
it has degree  less than $d$ for some $m>1$. 
This means, $\alpha_i^m=\alpha_j^m$ for some $1\leq i,j \leq d$. 
The lemma claims that the above defect can be cured by adding a suitable integer 
$k$ to $\alpha$. We will show that almost all integers $k$ will do. 
If a specific $k$ is unable to cure this defect it follows that
for suitable roots of unity $\zeta_{ij},\eta_{ij}$ we have simultaneously 

\begin{align}
\alpha_i &=\zeta_{ij}\alpha_j \label{1}\\
(\alpha_i+k)& =\eta_{ij}(\alpha_j+k) \label{2}
\end{align}

Substituting \eqref{1} in  \eqref{2} we get 
\begin{align*}
\zeta_{ij} \alpha_j +k=\eta_{ij}(\alpha_j+k)
\end{align*}
Rewriting the above, for this specific $k$, we get
\begin{align}
(\zeta_{ij}-\eta_{ij})\alpha_j=(\eta_{ij}-1)k\label{3}
\end{align}
We claim that Equation \eqref{3} will be true only for a finite number of $k$'s 
thereby proving the theorem. 
Suppose not. Then  we will have infinitely many choices of  $k$ which will
make the RHS of \eqref{3} grow unboundedly with such $k$.
However we can easily see that the LHS will be bounded with the following argument:
For all the choices of $k$ in \eqref{3} 
LHS will change only  in the choice of the roots of unity and the choice
among the $d$ conjugates of $\alpha$.
So  applying triangle inequality to  the LHS we get,
\[\abs{(\zeta_{ij}-\eta_{ij})\alpha_j} \leq  2\max \big\{|\alpha_1|,|\alpha_2|,\ldots,|\alpha_d|\big\}\]
which is an absolute bound independent of $k$ depending only on $\alpha$. \hfill QED

\vskip6pt

The following elementary fact about linear disjointness plays a very
significant role in our paper.

\vskip6pt
\noindent
 \textbf{Lemma 2:}\quad (See Morandi \cite{morandi}, page 184, Example 20.6)
 \textit{Suppose that $K$ and $L$  are finite extensions of a field $F$ with $K$ a Galois
 extension over $F$.
Under this hypothesis,  $K$ and $L$ are linearly disjoint  over $F$ 
if and only if $K\cap L=F$.} 

\vskip6pt

The following lemma which must be well-known to experts
is provided with a proof for the  reason that it is short, and also it will
make our  main theorem widely applicable.

\vskip6pt\noindent
\textbf{Lemma 3}: \textit{There exists number fields $F$ 
such that $F$ is linearly disjoint  over $\bq$ with every cyclotomic extension.}

\vskip6pt
Proof: 
We first claim that any number field $F$ which is a Galois extension of $\bq$
having some nonabelian simple group $G$ as Galois group will do. 
If such an $F$ is not linearly
disjoint with a cyclotomic extension $L$ then, by Lemma 2 it follows
 that $L\cap F=K$  will be a number field of degree${}>1$. Being contained
 in a cyclotomic extension, $K$ will be a Galois extension of $\bq$.
 As it is also a subfield of $F$,  by fundamental Galois correspondence,
 this $K$ must be a fixed field of a \textit{normal} subgroup of $G$.
 The simplicity of $G$ rules this out providing  the contradiction.

 To finish the argument we need to show that such an extension 
 $F$ with a simple non-abelian Galois group exists. Schur (see e.g.\cite{banerjee})
 has shown that for some  generalized Laguerre
 polynomials $L_n^{(\alpha)}(x),\,\alpha =1$ their splitting fields have 
 the alternating group $A_n$  as Galois groups. 
 Also, the specific polynomial $x^7-154x+99$ has been shown to have
 the simple group $PSL(2,\mathbf{F}_7)$ as the Galois group
  by Erbach, Fisher and McKay in \cite{efm}.

  \vskip4pt
  It seems there is no name for an algebraic number which generates
  a Galois extension of $\bq$.  Nor does the minimal  polynomial for such a number
  have any special name.

  As we will be constructing such polynomials later,
  we honour them formally here with a name:

\vskip6pt
\noindent
 \textbf{Galois polynomial:} 
 An irreducible polynomial $f(x)\in \bq [x]$  with the property that
  for any  root $\alpha$ of $f$,
 the field $\bq[\alpha] $ is also a splitting field of $f(x)$ over $\mathbf{Q}$ 
 is called a  \textit{Galois polynomial}.

Evidently,  cyclotomic polynomials are  Galois polynomials.
\section{Relative Cyclotomic Polynomials}

We state our main theorem under milder assumptions, though one can
show lots of examples where stronger hypothesis (and consequently stronger conclusions)
are valid.

\vskip6pt
\noindent
\textbf{Main Theorem}:\quad
\textit{Let $\alpha$ be an algebraic integer of degree $k$ and $f(x)$  its minimal polynomial
over $\bq$.
Choose  an integer $m > 2$ such that
\begin{itemize}
	\item[(i)] The $m$-th cyclotomic field $\mathbf{Q}[\zeta]$  and $\mathbf{Q}[\alpha]$ are linearly disjoint over
		$\bq$. (here $\zeta = \exp(2\pi i/m)$).
	\item[(ii)] the degree of $\alpha ^m$ is the same as the degree of $\alpha$.
\end{itemize}
 Denoting by  $f_m(x)$  the minimal polynomial of $\alpha^m$  over $\bq$,
  define $F_m(x) \mathrel{:=}f_m(x^m)$. 
 With this notation, 
  $F_m(x)$ factories in $\bq[x]$ with exactly one irreducible factor, denoted $\Phi_{d,f}(x)$,
  with  degree${}=\deg f\varphi(d) $   for each divisor $d$ of $m$. That is,}
 \[F_m(x)=f_m(x^m)=\prod\limits_{d|m} \Phi_{d,f}(x)\quad  \mbox{where}\quad 
 \deg(\Phi_{d,f})=\deg(f).\varphi(d).\]

\textrm{We call these polynomials $\Phi_{d,f}(x)$ as relative $d$-th cyclotomic polynomials, relative to 
 $f(x)$.}

\vskip6pt
 \textbf{Proof:}\quad
 Let $\alpha_1=\alpha,\alpha_2,\ldots,\alpha_k$ be the conjugates of $\alpha$.
  By definition, $F_m(x) $ has degree $km$. 
 The  $k$ conjugates of $\alpha ^m$, namely  $\alpha_1^m,\alpha_2^m,\ldots,\alpha_k^m$  
are  among the  $km$  roots of  $F_m(x)$. 
We first claim that the roots of $F_m(x)$ are
all distinct. In fact for any irreducible polynomial $h(x)$ 
in characteristic zero  $h(x^m)$ has distinct roots (because $h$ has).

Let $S$ denote the set of all those $km$ roots of $F_m$. Clearly
\[S=\{\alpha_i\zeta^j \mid i=1,2,\ldots,k,\ \mbox{and }j=1,2,\ldots,m \}.\]
Our plan is to partition $S$ into disjoint subsets $S_d$, one for each
divisor $d$ of $m$, with $S_d$ a  set of cardinality $k\varphi(d)$  
as the roots of an irreducible polynomial $\Phi_{d,f}(x)$.

With this in mind, for a $d\mid m$ we define 
\[S_d=\big\{\alpha_i\zeta_d^j \mid i=1,2,\ldots,k\  \mbox{and }
\gcd(j,d)=1\  \mbox{with} \ j<d \big \}\]

Then $S=\bigcup_{d\mid m} S_{d}$ . We can see that for any two distinct
divisors $d,d'$ of $m$ we have $S_d \cap S_{d'}=\emptyset $ and $\abs{S_d}=k\varphi(d)$.
 
 \[| S |  =\sum \limits _{d\mid m} |S_d|  = \sum \limits _{d\mid m} k\varphi(d)=
 k\sum \limits _{d\mid m} \varphi(d)=km.\]

 We denote by $\zeta_d$ a primitive $d$-th root of unity, except for the 
 case $\zeta_m$ we use simply $\zeta$.

Now we proceed to show that elements of $S_d$ are the roots of
the minimal  polynomial of $\alpha\zeta_d$.

Given $\bq[\alpha]$ and $\bq[\zeta]$ are linearly disjoint over $\bq$. 
Consequently, for any $d\mid m$, we have  $\bq[\alpha]$ and $\bq[\zeta_d]$ 
are also linearly disjoint.
So, we have the following equality of field degrees:
\[\big[\bq[\alpha,\zeta_d] :\bq \big]=\big[\bq[\alpha]:\bq \big]\times \big[\bq[\zeta_d]:Q \big]=k\varphi(d).\]
We know that $\bq[\alpha \zeta _d]\subset \mathbf{Q}[\alpha,\zeta_d]$.
So, $\big[\bq[\alpha \zeta_d] :\bq \big]\leq \big[\bq[\alpha,\zeta_d] :\mathbf{Q}\big]= k\varphi(d)$. 
If we show $\alpha \zeta_d$ is an algebraic number of degree $k\varphi(d)$ we are done.
This is equivalent to showing $\alpha\zeta_d$ has relative degree $\varphi(d)$ over
$\bq[\alpha]$. Suppose, to the contrary, this relative degree $r
 < \varphi(d)$ over $\mathbf{Q}[\alpha]$.
 Then 
 for some  $h(x)=\sum \limits _{i=0}^r c_i x^i$ with $c_i \in \mathbf{Q}[\alpha]$, we  have 
 $h(\alpha \zeta_d)=0$. That is, 
 \[0=\sum \limits _{i=0} ^r c_i(\alpha \zeta_d)^i=\sum \limits _{i=0} ^r
 (c_i\alpha^i) \zeta_d^i=\sum \limits _{i=0} ^r b_i \zeta_d^i\mbox{ for}\ b_i \in \bq[\alpha].\]
 This is a contradiction to the linear disjointness of $\bq[\alpha]$ with
 $\bq[\zeta_d]$.
  So $\alpha \zeta_d $ is an algebraic number of degree $\varphi(d)$ over $\mathbf{Q}[\alpha]$. Its  minimal polynomial over $\bq$ is the polynomial denoted by  $\Phi_{d,f}(x)$ 
  in the statement of the theorem.
  This is irreducible and is of degree $k\varphi(d)$.\hfill QED

\textbf{Corollary:} If, further, $\mathbf{Q}[\alpha]$ over $\mathbf{Q}$ is a Galois extension, then for each divisor $d|m$, the above polynomial  $\Phi_{d,f}(x)$ will be  a Galois polynomial
of degree $\varphi(d)\deg (f) $.

Now we would like to point out the scenarios where the hypothesis of our main theorem hold.

We will choose $\alpha $ to be of uniform degree, by adding an integer
(provided by Lemma 1)  so that condition (ii) of Main theorem
will be true \textit{for all}  $m$. To find examples where condition (i) holds
we can use examples of fields assured by  Lemma 3.

But for cases where the splitting field of $f(x)$ intersects with some 
cyclotomic extension, one can still find infinitely many
$m$ such that $m$-th cyclotomic field to be linearly disjoint with it.

For example, if we choose $\alpha = 1+\sqrt2$, the factorization
fails (that is, $F_m(x)$ does not imitate the factorization behaviour
of $x^m-1$ ) whenever $m$ is
a multiple of the discriminant (8 in this case).  Otherwise the factorization
will be along the classical lines.

\noindent\textbf{Remark}:\quad The relative cyclotomic polynomials $\Phi_{d,f}(x)$
relative to $f(x)$ (using the same notations as in the Main Theorem)
may alternatively be expressed using M\"obius inversion by
\[ \Phi_{m,f}(x) = \prod_{ d|m} F_d(x)^{\mu (n/d)} 
 =\prod_{d|m} f_d(x^d)^{\mu(n/d)}\]
\section{Examples and Applications}
Using our main theorem, given a Galois polynomial we can indicate how to find more 
Galois polynomials of higher degrees. More specifically,  starting with
Galois polynomial $f$ of degree $k$ with Galois group $G$ of order $k$, we can
	show how to find Galois polynomials 
	 Galois group of the form $G\oplus H$.
	 
 First requirement is  $H$ must be abelian, and it should be
 the Galois group of a suitable cyclotomic extension of $\bq$
 
\textit{ All the calculations were done with the computer algebra system SAGE.}

\vskip6pt
\noindent  \textbf{EXAMPLE 1}

 We start with  $f=x^6 + 3x^5 - 2x^4 - 9x^3 + 5x + 1$ a Galois polynomial 
 with Galois group $S_3$ found in \textit{A Database for Number Fields} in the web site
	\texttt{http://galoisdb.math.upb.de/home} as the base polynomial.

We illustrate this first for $H=C_2\oplus C_2$, the Klein's group
which is the Galois group of the 8-th cyclotomic extension.
	As this extension is linearly disjoint with $\bq[\alpha]$
 ($\alpha$  denoting a root of this $f$), we can appeal to our Main Theorem.
 So we look at $\alpha^8$ which is again an algebraic number of degree $6$ with minimal polynomial $f_8(x)$. From the factorization of $F_8(x)=f_8(x^8)$ we obtain new Galois polynomial, namely the unique relative cyclotomic polynomial 
 factor $\Phi_{8,f}(x)$:
 
\begin{align*} F_8(x)=(x^{48} - 1405x^{40} + 226310x^{32} - 3670777x^{24}
 + 1940230x^{16}-201085x^8 +1) 
\end{align*}

\begin{align*}
	\Phi_{8, f}(x)=  (x^{24} + 53x^{20} + 702x^{16} + 2553x^{12} + 2062x^8 + 453x^4 + 1)  
 \end{align*}
The polynomial $\Phi_{8,f}(x)$ is of degree $24=\deg f\cdot\varphi(8)$ and  has Galois group 
$\mathrm{Gal\,}\big(f(x)\big)\oplus \mathrm{Gal\,}\big(\Phi_8(x)\big)$
of order 24 which is isomorphic to $S_3\oplus  C_2 \oplus C_2$.
This is a Galois polynomial.

\vskip6pt 
\noindent  \textbf{EXAMPLE 2}

  For the Galois polynomial $g=(x^6+3x^5+5x^4+5x^3+5x^2+3x+1)$ with Galois group $S_3$,
  found in the same database
 we can compute and arrive at  \[F_3(x)= (x^{18} - 3x^{15} + 11x^{12} + 5x^9 + 11x^6 - 3x^3 + 1)\] 
which leads us to the following relative cyclotomic polynomial as its highest degree 
factor.
\begin{align*}
\begin{aligned}
	\Phi_{3,g}(x) &= (x^{12} - 3x^{11} + 4x^{10} - 5x^9 + 5x^8 + 2x^7 - 7x^6 + 2x^5 + 5x^4\\&-5x^3 +4x^2 - 3x + 1)
\end{aligned}
\end{align*}

 The relative cyclotomic factor $\Phi_{3,g} (x)$ is of degree 
 $12 = \deg g\cdot \varphi(3)$ and 
 has Galois group $\mathrm{Gal\,}\big(g(x)\big)\oplus \mathrm{Gal\,}\big(\Phi_3(x)\big)$
 of order 12 isomorphic to $S_3\oplus C_2$. This is a Galois polynomial.

\vskip6pt 
\noindent  \textbf{EXAMPLE 3}

Using the same $g$ as in Example 2, but with $m=5$ this time, we obtain

  $F_5(x)=(x^{30} - 17x^{25} + 95x^{20} - 135x^{15} + 95x^{10} - 17x^5 + 1)$\\
\begin{align*}
 \Phi_{5,g}(x)&= (x^{24} - 3x^{23} + 4x^{22} - 2x^{21} - 4x^{20} - 3x^{19} + 27x^{18} - 45x^{17}\\
  &+ 37x^{16} + 8x^{15} - 11x^{14} - 45x^{13} + 73x^{12} - 45x^{11} - 11x^{10} + 8x^9 + 37x^8 \\
  &- 45x^7 + 27x^6 - 3x^5 - 4x^4 - 2x^3 + 4x^2 - 3x + 1)
 \end{align*}
 The relative cyclotomic factor $\Phi_{5,g}(x)$ is a polynomial
 of degree $24=\deg g\cdot\varphi(5)$, has Galois group $\mathrm{Gal\,} \big(g(x)\big)\oplus
 \mathrm{Gal\,}\big(\Phi_5(x)\big)$ of order  24 which is  isomorphic to $S_3\oplus C_4$.

\section{Conclusion}
When the linear disjointness hypothesis with the $m$-th cyclotomic field
does not hold, it was found empirically that
we get excess factors for $F_m(x)$.
To elaborate,  in all the examples computed for those
failed cases we found that the number of irreducible
factors of $F_m(x)$ was
\textit{always} more than the number of divisors of $m$.
   This we are yet to prove and needs further  investigation.

\end{document}